\documentclass[leqno]{beamer}
\usepackage{newlattice} 

\DeclareMathOperator{\Princ}{Princ}
\DeclareMathOperator{\Conc}{\tup{Con}_\tup{c}}
\newcommand{\Pd}{P^{\,\tup{d}}}

\begin{document} 
\title[The order of principal congruences of a bounded lattice]{The order of principal congruences\\ of a bounded lattice.\\
AMS Fall Southeastern Sectional Meeting\\
University of Louisville, Louisville, KY\\ 
October 5-6, 2013}
\author{G. Gr\"atzer}
\date{}
\maketitle 

\begin{frame}
\frametitle{Summary}
We characterize the order of principal congruences 
of a bounded lattice 
as a bounded ordered set.
We~also state a number of open problems in this new field.
\medskip

arxiv: 1309.6712 
\end{frame}

\begin{frame}
\frametitle{Background}
Let $A$ be a lattice (resp., join-semilattice with zero). 
We call $A$ \emph{representable} 
if there exist a lattice $L$
such that $A$ is isomorphic to the congruence lattice of $L$, 
in formula, $A \iso \Con L$ 
(resp., $A$ is isomorphic to the join-semilattice with zero
of compact congruences of $L$, 
in~formula, $A \iso \Conc L$).
\end{frame}

\begin{frame}
\frametitle{Background}
For over 60 years, one of lattice theory's most central conjectures 
was the following:
\begin{quote}
\emph{Characterize representable lattices as distributive algebraic lattices.}
\end{quote}
\pause
Or equivalently:
Characterize representable join-semilattices as distributive join-semilattice with zero.
\pause

This conjecture was refuted in F. Wehrung in 2007.
\end{frame}

\begin{frame}
\frametitle{Principal congruences}
In this lecture, we deal with $\Princ L$, 
the order of principal congruences of a lattice $L$.
Observe that 
\begin{enumeratea}

\item $\Princ L$ is a directed order with zero.
\pause
\item $\Conc L$ is the set of compact elements of $\Con L$, 
a lattice theoretic characterization of this subset.
\pause
\item $\Princ L$ is a directed subset of $\Conc L$ 
containing the zero and join-generat\-ing $\Conc L$;
there is no lattice theoretic characterization of this subset.
\end{enumeratea}
\end{frame}

\begin{frame}
\frametitle{Principal congruences}
\centerline{\includegraphics[scale=.6]{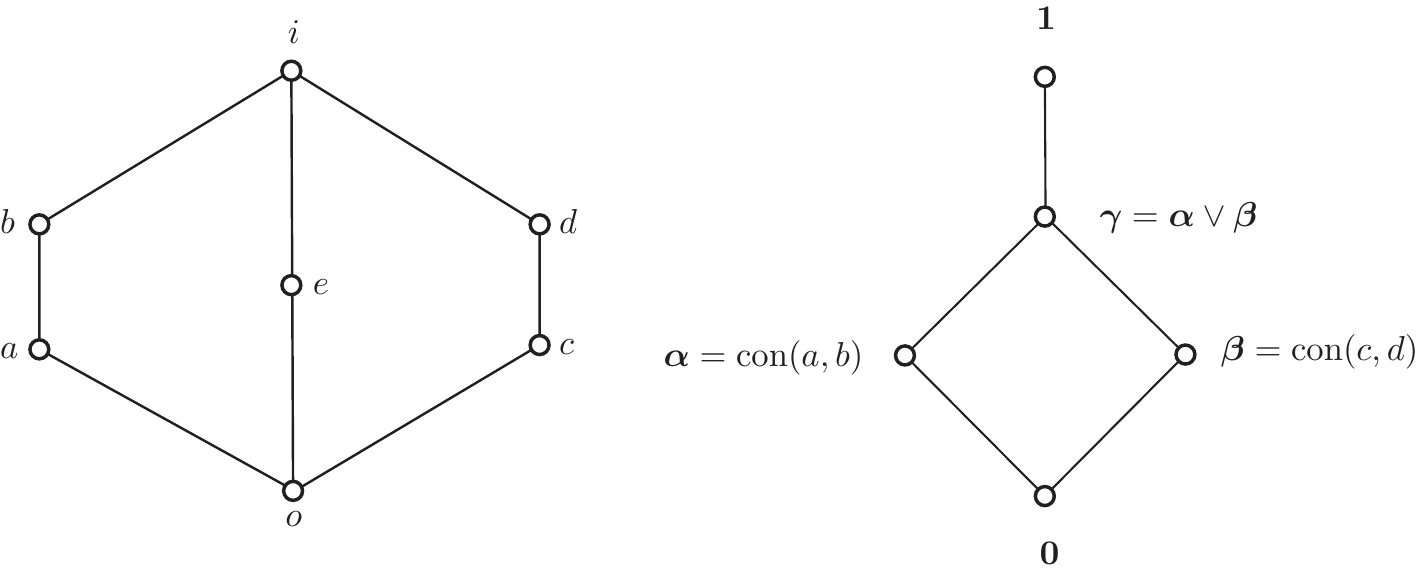}}

This is the lattice $\SN 7$ 
and its congruence lattice $\SB 2 + 1$.

Note that $\Princ \SN 7 = {\Con \SN 7} - \set{\bgg}$,
while in the standard representation~$K$
of $\SB 2 + 1$ as a congruence lattice
(G. Gr\"atzer and E.\,T. Schmidt, 1962), 
we have $\Princ K = \Con K$.
 
 \pause
This example shows that $\Princ L$
has no lattice theoretic description in $\Con L$.
\end{frame}

\begin{frame}
\frametitle{Theorem 1}
For a bounded lattice $L$, the order $\Princ K$ is bounded. 
We now state the converse.
\pause
\begin{theorem}
Let $P$ be an order with zero and unit.
Then there is a bounded lattice~$K$ such that
\[
   P \iso \Princ K.
\]
If $P$ is finite, we can construct $K$ as a finite lattice.
\end{theorem}

\end{frame}

\begin{frame}
\frametitle{Lattice Problem 1}
\begin{problem}
Can we characterize the order $\Princ L$ for a lattice $L$
as a directed order with zero? 
\end{problem}

\pause
(Why directed?)
\pause

G. Cz\'edli solved this problem for countable lattices 

arXiv:1305.0965
\end{frame}

\begin{frame}
\frametitle{Lattice Problem 2}
Even more interesting would be to charaterize 
the pair $P = \Princ L$ in $S =\Conc L$ by the properties
that $P$ is a directed order with zero 
that join-generates~$S$. 
We have to rephrase this 
so it does not require a solution 
of the congruence lattice characterization problem.
\pause
\begin{problem}
Let $S$ be a representable join-semilattice. 
Let $P \ci S$ be a directed order with zero 
and let $P$ join-generate~$S$.
Under what conditions is there a lattice $K$ such 
that $\Conc K$ is isomorphic to $S$
and under this isomorphism $\Princ K$ corresponds to $P$?
\end{problem}
\end{frame}

\begin{frame}
\frametitle{Lattice Problem 3}
For a lattice $L$, let us define a \emph{valuation} $v$ 
on $\Conc L$ as follows: %
for a compact congruence \bga of $L$, 
let $v(\bga)$ be the smallest integer $n$
such that the congruence \bga is the join of $n$ principal congruences.

A valuation $v$ has some obvious properties, for instance,
$v(\zero) = 0$ and $v(\bga \jj \bgb) \leq v(\bga) + v(\bgb)$.
Note the connection with $\Princ L$:
\[
   \Princ L = \setm{\bga \in \Conc L}{v(\bga) \leq 1}.
\]
\pause
\begin{problem}
Let $S$ be a representable join-semilattice. 
Let $v$ map $S$ to the natural numbers.
Under what conditions is there an isomorphism \gf
of $S$ with $\Conc K$ for some lattice $K$
so that under \gf the map $v$ corresponds 
to the valuation on $\Conc K$?
\end{problem}
\end{frame}

\begin{frame}
\frametitle{Lattice Problem 4}
Let $D$ be a finite distributive lattice.
In G. Gr\"atzer and E.\,T. Schmidt 1962,
we represent $D$ as the congruence lattice
of a finite lattice~$K$ in which 
\emph{all congruences are principal} (that is, $\Con K = \Princ K$).
\pause
\begin{problem}
Let $D$ be a finite distributive lattice. 
Let $Q$ be a subset of $D$ satisfying
$\set{0, 1} \uu \Ji D \ci Q \ci D$.
When is there a finite lattice~$K$ such 
that $\Con K$ is isomorphic to $D$
and under this isomorphism $\Princ K$ 
corresponds to $Q$?
\end{problem}
\end{frame}

\begin{frame}
\frametitle{Lattice Problem 4, an example}

Example:

Let $D$ be the eight-element Boolean lattice.
Let $Q$ be a subset of~$D$ containing
$0$ and $1$ and the three atoms (the join-irreducible elements).

\begin{lemma}
If there is a finite lattice~$K$ such 
that $\Con K$ is isomorphic to $D$
and under this isomorphism $\Princ K$ 
corresponds to $Q$,
then $Q$ has seven or eight elements.\end{lemma}

\end{frame}

\begin{frame}
\frametitle{Lattice Problem 5}
In particular, let $Q = \Con L$.

\begin{problem}
Let $\bold K$ be a class of lattices
with the property that every finite distributive lattice $D$
can be represented 
as the congruence lattice of some finite lattice in $\bold K$. 
Under what conditions on $\bold K$ is it true that every 
every finite distributive lattice $D$
can be represented 
as the congruence lattice of some finite lattice $L$ in $\bold K$
with the additional property: $\Con L = \Princ L$.
\end{problem}
\end{frame}

\begin{frame}
\frametitle{Theorem 2}

G. Gr\"atzer and E.\,T. Schmidt, 
\emph{An extension theorem for planar semimodular lattices.}
Periodica Mathematica Hungarica.
arXiv: 1304.7489

\begin{theorem}
Every finite distributive lattice $D$
can be represented as the congruence lattice of 
a finite, planar, semimodular lattice $K$ with the property that
all congruences are principal.
\end{theorem}

In fact, $K$ is constructed as a ``rectangular lattice''.
\end{frame}

\begin{frame}
\frametitle{Problem 6}
In the finite variant of the valuation problem,
we need an additional property.

\begin{problem}
Let $S$ be a finite distributive lattice. 
Let $v$ be a map of $D$ to the natural numbers
satisfying $v(0) = 0$, $v(1) = 1$, and 
$v(a \jj b) \leq v(a) + v(b)$ for $a, b \in D$.
When is there an isomorphism \gf
of $D$ with $\Con K$ for some finite lattice~$K$
such that under \gf the map $v$ corresponds 
to the valuation on $\Con K$?
\end{problem}
\end{frame}

\begin{frame}
\frametitle{Problem 7}

\begin{problem}
In Theorem 1, can we construct a semimodular lattice?
\end{problem}

\pause

Remember Theorem 1:

\begin{theorem}
Let $P$ be an order with zero and unit.
Then there is a bounded lattice~$K$ such that
\[
   P \iso \Princ K.
\]
If $P$ is finite, we can construct $K$ as a finite lattice.
\end{theorem}
\end{frame}

\begin{frame}
\frametitle{Problem 8}
\begin{problem}
In Problems 2 and 3, in the finite case, 
can we construct a finite semimodular lattice $K$?
\end{problem}

\pause

Remember Problems 2 and 3:

\begin{problem}
Let $S$ be a representable join-semilattice. 
Let $P \ci S$ be a directed order with zero 
and let $P$ join-generate~$S$.
Under what conditions is there a lattice $K$ such 
that $\Conc K$ is isomorphic to $S$
and under this isomorphism $\Princ K$ corresponds to $P$?
\end{problem}

\begin{problem}
Let $S$ be a representable join-semilattice. 
Let $v$ map $S$ to the natural numbers.
Under what conditions is there an isomorphism \gf
of $S$ with $\Conc K$ for some lattice $K$
so that under \gf the map $v$ corresponds 
to the valuation on $\Conc K$?
\end{problem}
\end{frame}

\begin{frame}
\frametitle{Problem 9}
In E.~T. Schmidt 1962
(see also G. Gr\"atzer and E.\,T. Schmidt 2003),
for a finite distributive lattice $D$, 
a countable modular lattice $M$ is constructed
with $\Con M \iso D$.

\begin{problem}
In Theorem 1, for a finite $P$,
can we construct a countable modular lattice $K$?
\end{problem}
\end{frame}

\begin{frame}
\frametitle{}
Some of these problems seem to be of interest for algebras
other than lattices as well.

\begin{problem}
Can we characterize the order $\Princ \F A$ 
for an algebra~$\F A$
as an order with zero?
\end{problem}
\end{frame}

\begin{frame}
\frametitle{Problem 9}

\begin{problem}
For an algebra $\F A$,
how is the assumption that the unit congruence~$\one$
is compact reflected in  the order $\Princ \F A$? 
\end{problem}
\end{frame}

\begin{frame}
\frametitle{Problem 10}

\begin{problem}
Let  $\F A$ be an algebra
and let $\Princ \F A \ci Q \ci \Conc \F A$. 
Does there exist an algebra $\F B$ such that
$\Con \F A \iso \Con \F B$ and under this isomorphism
$Q$ corresponds to $\Princ \F B$? 
\end{problem}
\end{frame}

\begin{frame}
\frametitle{Problem 11}

\begin{problem}
Extend the concept of valuation to algebras in general.
State and solve Problem~3 for algebras.
\end{problem}

\pause

Remember Problem 3:

\begin{problem}
Let $S$ be a representable join-semilattice. 
Let $v$ map $S$ to the natural numbers.
Under what conditions is there an isomorphism \gf
of $S$ with $\Conc K$ for some lattice $K$
so that under \gf the map $v$ corresponds 
to the valuation on $\Conc K$?
\end{problem}
\end{frame}

\begin{frame}
\frametitle{Problem 12}

\begin{problem}
Can we sharpen the result of 
G. Gr\"atzer and E.\,T. Schmidt 1960:
every algebra $\F A$ 
has a congruence-preserving extension $\F B$
such that $\Con \F A \iso \Con \F B$  
and $\Princ \F B = \Conc {\F B}$.
\end{problem}
\pause
I do not even know whether every algebra $\F A$ has a 
proper congruence-preserving extension $\F B$.
\end{frame}

\begin{frame}
\frametitle{Proof by Picture}
For a bounded order $Q$, 
let $Q^-$ denote the order $Q$ with the bounds removed.
Let $P$ be the order in Theorem~1. 
Let $0$ and $1$ denote the zero and unit of $P$, respectively.
We denote by $\Pd$ those elements of $P^-$ that are not comparable to any other element of $P^-$, that is,
\[
   \Pd = \setm{x \in P^-}{ x \parallel y 
      \text{ for all } y \in P^-,\ y \neq x}.
\]

\end{frame}
\begin{frame}
\frametitle{Proof by Picture: The Lattice $F$}
We first construct the lattice $F$
consisting of the elements \text{$o$, $i$} 
and the elements $a_p, b_p$ for every $p \in P$,
where $a_p \neq b_p$ for every $p \in P^-$
and $a_0 = b_0$, $a_1 = b_1$.  

The lattice $F$:
\begin{figure}[h!]
\centerline{\includegraphics[scale=0.9]{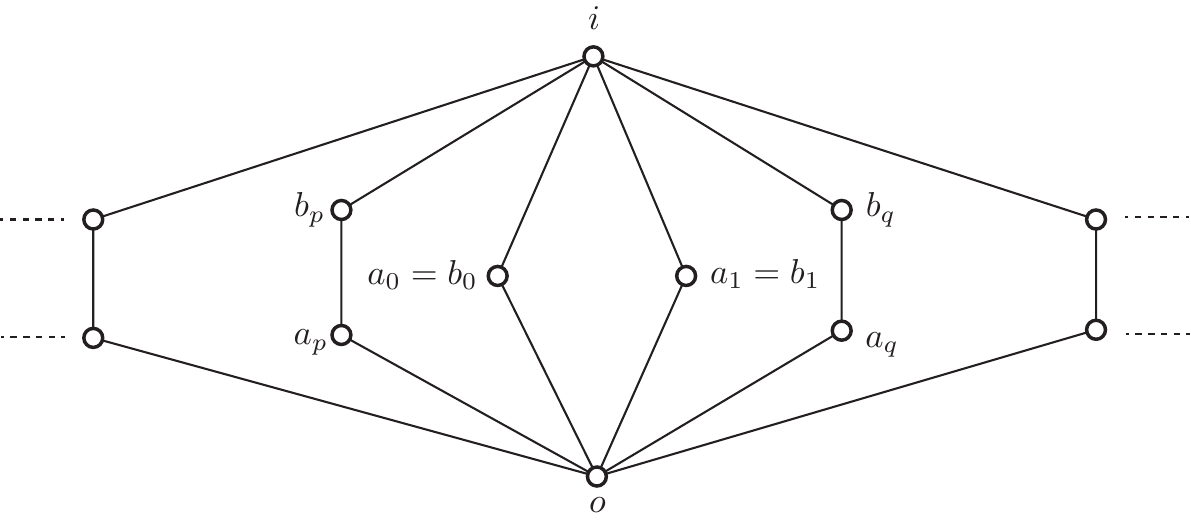}}
\end{figure}

\end{frame}

\begin{frame}
\frametitle{Proof by Picture: The Lattice $K$}
We are going to construct the lattice $K$ (of Theorem~1)
as an extension of $F$. 
For $p \prec q$, between the edges $[a_p, b_p]$ and $[a_q, b_q]$
we insert the lattice $S = S(p, q)$:

\begin{figure}[hbt]
\centerline{\includegraphics[scale=0.7]{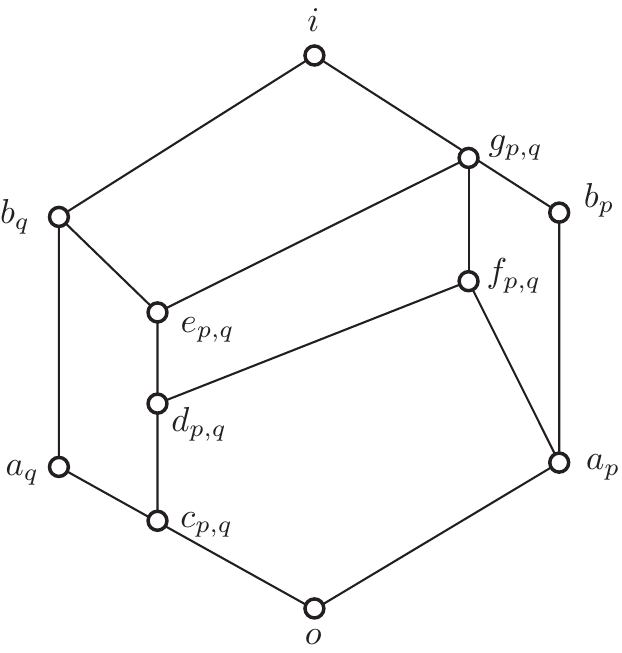}}\label{F:S}
\end{figure}

The principal congruence of $K$ representing $p \in P^-$ 
will be $\con{a_p, b_p}$.

\end{frame}

\begin{frame}
\frametitle{Proof by Picture: The Orders $C$, $V$, and $H$}
For $x \in S(p, q)$ and $y \in S(p', q')$, $p \prec q$,  $p' \prec q'$
we have to find
$x \jj y$ and $x \mm y$. 

If $\set{p, q} \ii \set{p', q'} = \es$, then $x$ and $y$
are complimentary.

If $\set{p, q} \ii \set{p', q'} \neq \es$, then 
$\set{p, q} \uu \set{p', q'}$ form a three element order
$\SC{}$, $\SV{}$, or $\SH{}$:

\begin{figure}[p]
\centerline{\includegraphics[scale=.85]{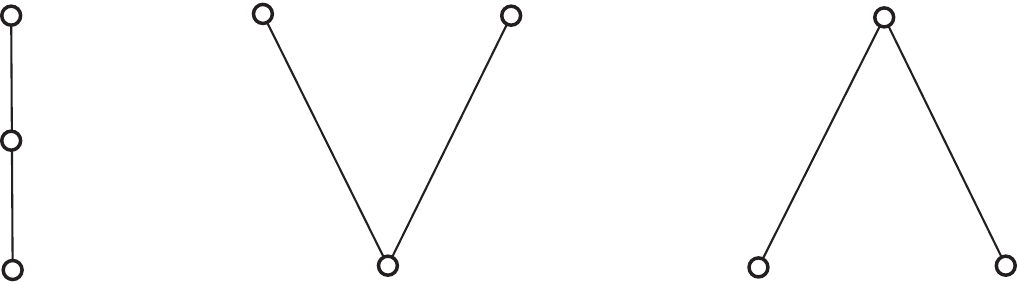}}
\end{figure}
\end{frame}

\begin{frame}
\frametitle{Proof by Picture}

We form $x \jj y$ and $x \mm y$ in the appropriate lattices, 
$S_{\SC{}} = S(p < q,\ q < q')$,
$S_{\SV{}} = S(p < q,\ p < q')$ with $q \neq q'$, and
$S_{\SH{}} = S(p < q,\ p' < q)$ with $p \neq p'$.
\end{frame}

\begin{frame}
\frametitle{Proof by Picture: The Lattice $S_{\SC{}}$}
The lattice $S_{\SC{}} = S(p < q,\ q < q')$:
\begin{figure}[p]
\centerline{\includegraphics[scale=.85]{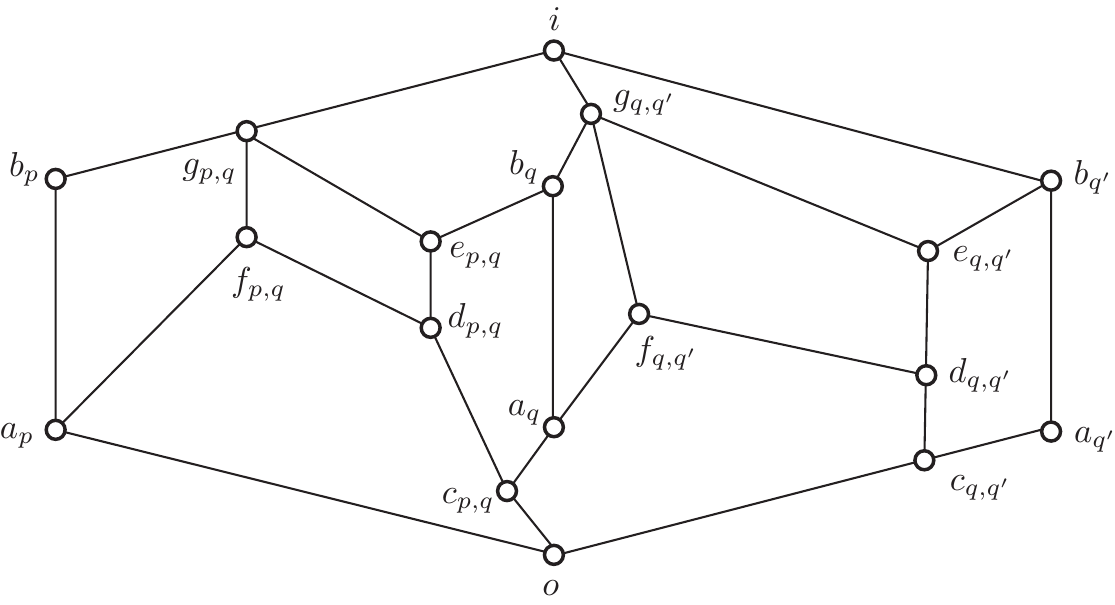}}
\end{figure}

\end{frame}

\begin{frame}
\frametitle{Proof by Picture: The Lattice $S_{\SV{}}$}
The lattice $S_{\SV{}} = S(p < q,\ p < q')$ with $q \neq q'$:
\begin{figure}
\centerline{\includegraphics[scale=.85]{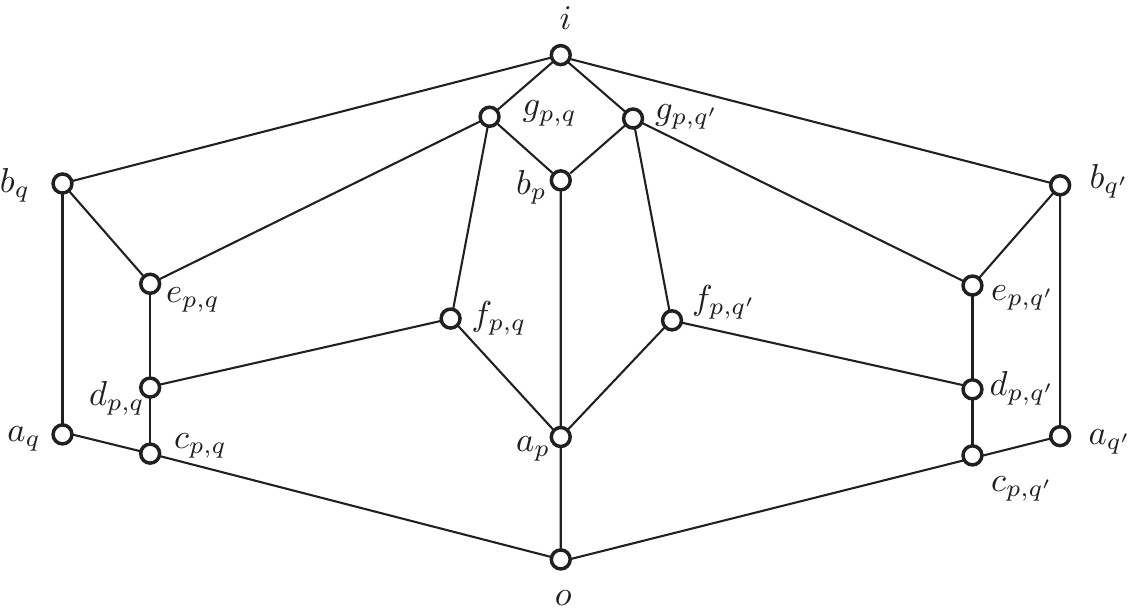}}
\end{figure}

\end{frame}

\begin{frame}
\frametitle{Proof by Picture: The Lattice $S_{\SH{}}$}
The lattice $S_{\SH{}} = S(p < q,\ p' < q)$ with $p \neq p'$:
\begin{figure}
\centerline{\includegraphics[scale=.85]{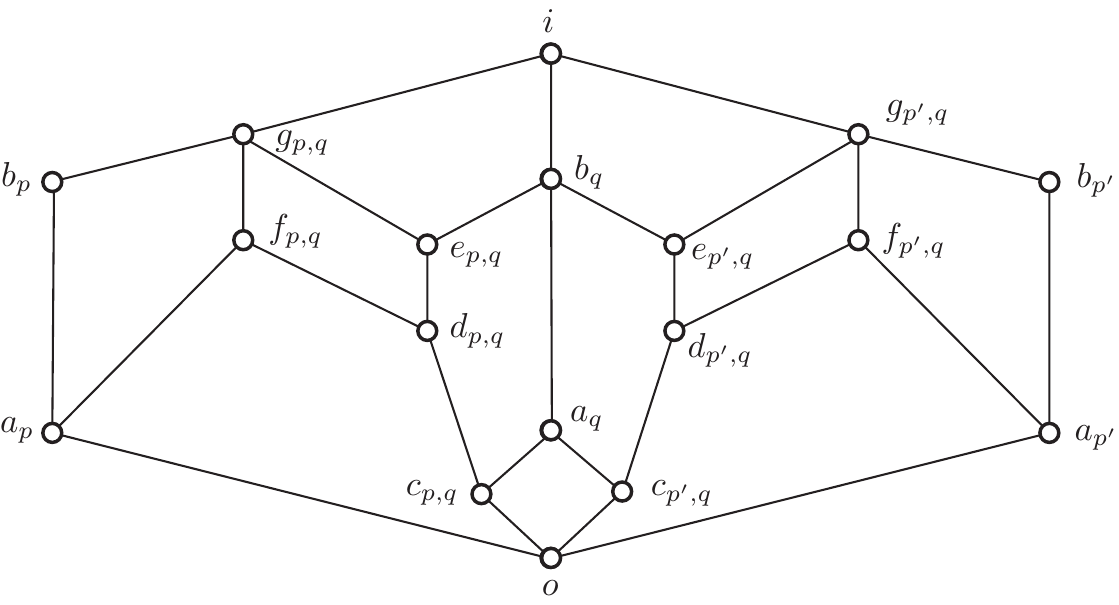}}
\end{figure}

\end{frame}

\end{document}